\documentclass{amsart}
\usepackage{amsfonts,amscd}
               
\newtheorem{claim}{}[section]
\newtheorem{theorem}[claim]{Theorem}
\newtheorem{lemma}[claim]{Lemma}
\newtheorem{proposition}[claim]{Proposition}
\newtheorem{corollary}[claim]{Corollary}

\newtheorem{definition}[claim]{Definition}

\def\proclaim #1. #2\par{\medbreak
\noindent{\bf#1.\enspace}{\sl#2}\par\medbreak}
\makeatother

\DeclareMathOperator{\Bdb}{{\mathbb B}}
\DeclareMathOperator{\Kdb}{{\mathbb K}} 
\DeclareMathOperator{\Cdb}{{\mathbb C}}

\DeclareMathOperator{\Ndb}{{\mathbb N}}

\DeclareMathOperator{\M}{{\mathcal M}}
\DeclareMathOperator{\A}{{\mathcal A}}

\begin{document}
                                        
\title[Multipliers between two operator spaces]{Multipliers between two operator spaces}
 
\date{August 23, 2003}
 
\author{David P. Blecher}
\address{Department of Mathematics, University of Houston, Houston,
TX
77204-3008}
\email[David P. Blecher]{dblecher@math.uh.edu}
\thanks{*This research was supported in part by a grant from
the National Science Foundation}    
 
\begin{abstract}
In a recent survey paper we introduced 
one-sided multipliers between two different 
operator spaces.   Here we  give some basic theory
for  these maps.      \end{abstract}
 
\maketitle

%%%%%%%%%%%%%%%%%%%%%%%%%%%%%%%%%%%%%%%%%%%%%%%%%%%%%%%%%%%%%%%%%%%%%%
 
\let\text=\mbox

\bigskip
 
\section{Introduction}
 
In \cite{BShi} we introduced 
the space $\M_l(X)$ (resp. $\A_l(X)$) of
 left multipliers (resp.\ adjointable multipliers)
on an operator space $X$.  See also
\cite{Wer}.   This is an operator algebra (resp.\  $C^*$-algebra)
whose elements are maps on $X$; namely
the maps satisfying variants of 
Theorems \ref{chmu2} or \ref{newest} (resp.\ \ref{genad})
 below (in the case that $X = Y$).  Theory and applications of 
these maps may be found in a series
of papers (e.g.\ \cite{BShi,BP,BEZ,Dual,BZ}).
A survey of this theory containing proofs and 
applications, may be found in \cite{Bsv}.
In the latter paper we also began to consider
one-sided multipliers between two different
operator spaces $X$ and $Y$.
In the present paper we give some basic theory of  these maps.

To have any interesting notion of multipliers between two different
spaces $X$
and $Y$ it does seem that there does need to be some relation between 
$X$ and $Y$,
or between two $C^*$-algebras that are
canonically associated with $X$ and $Y$.
This is the approach taken in Sections 3 and 4, where
we require the `noncommutative Shilov boundaries'
of $X$ and $Y$ to be a module over two equal, 
or comparable, $C^*$-algebras.   We may then prove 
characterizations of the associated
multipliers from $X$ to $Y$, which are 
satisfying except for the feature that they do 
depend on the particular Shilov boundaries chosen.
For this reason we call these {\em relative left multipliers}.

The problem in the last paragraph
also seems to be related to the issue of having a well defined
`column sum' of two operator spaces $X$ and $Y$.  
If $X, Y$ are subspaces of an operator space $V$
then this difficulty evaporates,
we may simply define the column sum $X \oplus_V Y$
to be the algebraic sum $X \oplus Y$ endowed with an operator
space structure by identifying it with a subspace of $C_2(V)$
via the map
$$(x,y) \mapsto \left[
\begin{array}{c}
x \\
y \end{array}
\right] , \qquad x \in X, y \in Y . $$
This is the approach taken in the final two sections of this
paper, and in this case we can remove 
the `dependence on the Shilov boundary' problem
mentioned above.

\section{Preliminaries}

We reserve the symbols $H, K$ for Hilbert spaces.
We use standard notation for 
operator spaces (see e.g.\ \cite{ERbook,Pau,Pis}).
  In particular we write $C_n(X)$ and $R_n(X)$
for $M_{n,1}(X)$ and $M_{1,n}(X)$ respectively.  
We recall that the space  $CB(X,Y)$ of completely bounded 
maps between two operator spaces is also an operator space.   
An {\em injective operator space} $X$ has the property that
any completely contractive $T : Y \to X$ has a 
completely contractive extension $\tilde{T} : 
Z \to X$.  Here $Z$ is any operator space with subspace $Y$.
Alternatively, a subspace $X \subset B(K,H)$ is injective if and 
only if there is a completely contractive projection of
$B(K,H)$ onto $X$.

We will need some basic theory of $C^*$-modules;
 see \cite{Lan} or Section 3 of \cite{Bsv}.
In fact we will need a few facts from the latter paper which 
we now repeat:
Suppose that $W$ and $Z$ are right $C^*$-modules over $C^*$-algebras
$B$ and $C$ respectively, where $B$ is a
$C^*$-subalgebra of $C$.
Then we may define an analogue of the
$C^*$-module direct sum.  Namely we define 
$W \oplus^C Z$ to be the algebraic sum
$W \oplus Z$, endowed with matrix norms 
$$\left| \left|  \left[ \begin{array}{cl} w_{ij} \\ z_{ij} \end{array}
\right] \right| \right|  \; = \;
\Big| \Big|  \Bigl[ \sum_{k=1}^n \langle w_{ki} | w_{kj}
\rangle + \langle z_{ki} | z_{kj}
\rangle \Bigr] \Big| \Big|^{\frac{1}{2}} , \qquad [w_{ij}] \in M_n(W), [z_{ij}]
 \in M_n(Z).
$$
There are many ways to see that $W \oplus^C Z$ is
an operator space: this is explained in \cite{Bsv} for example.
In fact this will be obvious in most of the cases we are interested in 
the later sections.

 The following is an extension of a result of Paschke
(see \cite[Theorem 2.8]{Pas}).
 
\begin{theorem} \label{PasB2}  \cite{Pas,Bsv}  Suppose that
$W$ and $Z$ are right $C^*$-modules over $C^*$-algebras
$B$ and $C$ respectively, where $B$ is a
$C^*$-subalgebra of $C$.
If  $u : W \to Z$ is a $\mathbb{C}$-linear map,
then the following are equivalent:
\begin{itemize}
\item [(i)]  $u$ is a contractive $B$-module map;
\item [(ii)]  $\langle u(w) | u(w) \rangle \leq
\langle w |  w \rangle$, for all $w \in W$;
\item [(iii)]  $\left| \left| \left[
\begin{array}{c}
u (w) \\
z
\end{array}
\right]  \right|  \right| \; \leq \;
 \left| \left| \left[
\begin{array}{c}
w \\
z
\end{array}
\right]  \right|  \right|$, for all $w \in W, z \in Z$.
\item [(iv)]   $\left| \left| \left[
\begin{array}{c}
u (w) \\
c
\end{array}
\right]  \right|  \right| \; \leq \;
 \left| \left| \left[
\begin{array}{c}
w \\
c
\end{array}
\right]  \right|  \right|$, for all $w \in W, c \in C$.
\end{itemize}
\end{theorem}

Suppose that $Y, Z$ are
$C^*$-modules over $B$, and that $J$ is a closed ideal of $B$
(or of the multiplier $C^*$-algebra $M(B)$ of $B$)
containing the span of the ranges of the $B$-valued
inner products on $Y$ and $Z$.  We write $B_B(Y,Z)$
for the set of bounded $B$-module maps from $Y$ to $Z$.
It clearly follows from Theorem \ref{PasB2},
although it is also easy to prove  directly using Cohen's 
factorization theorem, that $Y$ and $Z$ are
$C^*$-modules over $J$, and that
\begin{equation} \label{obre}
B_B(Y,Z) = B_J(Y,Z) .
\end{equation} 
  Thus there is not a truly essential dependence of 
$B_B(Y,Z)$ on $B$.  
 
We write $\Bdb(Y,Z)$  for the set of 
 {\em adjointable maps} $T : Y \to Z$, that is
those maps for which there is a map $S : Z \to Y$ such that
$$\langle Ty | z  \rangle \; = \; \langle y | S z  \rangle
, \qquad y \in Y, z \in Z.
$$                                             
The space $\Kdb(Y,Z)$, of so-called
`compact' adjointable maps, is the closure in $\Bdb(Y,Z)$
of the span of the operators of the form $y \mapsto z \langle w | y
\rangle$ on $Y$, for $w \in Y, z \in Z$.
The set $\Kdb(Y) = \Kdb(Y,Y)$
is a $C^*$-subalgebra of the $C^*$-algebra $\Bdb(Y)
= \Bdb(Y,Y)$.    

For any right $C^*$-module $Y$
over a  $C^*$-algebra $B$
we may define the {\em linking $C^*$-algebra} ${\mathcal L}(Y)$
to be $\Kdb(Y \oplus_c B)$.
See e.g.\ \cite[Section 3]{Bsv} for more details if needed.
Since $Y$ is a subspace of ${\mathcal L}(Y)$,
it is consequently an operator space.  
We will always suppose that any
 $C^*$-module met in this paper has this
 operator space structure.
  
A  $C^*$-module $W$ over a $W^*$-algebra $M$
such that $W$ has a predual Banach space will 
be called a {\em $W^*$-module}.  By a result of
Zettl (see  \cite{EOR} for a modern proof), this is exactly the 
(extremely important and well behaved)
class of $C^*$-modules over $M$ which are
{\em self-dual}.  See \cite{Pas} for the meaning 
of this term, and for basic theory
of such modules.

Due to the work of Arveson and Hamana
(see e.g.\ \cite{SOC,Ham2}) 
it is known that
every operator space $X$ has a {\em noncommutative
Shilov boundary}.  We will write this as $\partial X$ 
or $(\partial X,i)$, where $\partial X$ is
a right $C^*$-module and $i : X \to \partial X$
is a linear complete isometry.  We refer 
the reader to any one of \cite{Ham2,BShi,Bsv},
for example, for a better 
description of $\partial X$ and of its important
universal property.  We will not explicitly write down
the latter property here but, loosely speaking, this
universal property  says that
$\partial X$ 
is a {\em smallest} $C^*$-module
containing $X$ completely
isometrically.  To be symmetrical, one should perhaps say 
`$C^*$-bimodule' here instead of
$C^*$-module, but since we are emphasizing left multipliers
we shall think of $\partial X$ as a 
right $C^*$-module over a $C^*$-algebra $B$ say.
The $B$-valued inner product on $\partial X$ we shall write 
as $\langle \cdot | \cdot \rangle$; we refer to this
as a (right) {\em Shilov inner product}.  
The  ideal in $B$ densely spanned by the range of this
inner product will be written as ${\mathcal F}(X)$,
or ${\mathcal F}$ if $X$ is understood.
By the afore-mentioned 
universal property one may see that
the noncommutative Shilov boundary, the Shilov inner product,
and the $C^*$-algebra ${\mathcal F}$,
 are essentially unique up to an appropriate isomorphism. 

We review quickly Hamana's
method construction of a noncommutative Shilov boundary
for $X$. 
We begin with an {\em injective envelope} $(I(X),i)$ of 
$X$.  By this term, we mean that 
$I(X)$ is an injective operator space, that 
$i : X \to I(X)$ is a complete isometry, and that 
the identity map is the only completely contractive linear
map from $I(X)$ to itself extending the identity map on $X$.
This is called the {\em rigidity property} of the
injective envelope.  
In fact $I(X)$ may be chosen to be a full right $C^*$-module
over a $C^*$-algebra ${\mathcal D}$
(see e.g.\ \cite{Ham2};
also ${\mathcal D}$ is the algebra denoted $I(X)^* I(X)$ in
\cite[p.3]{BP}).     We let
$\partial X$ be the smallest closed 
subspace $E$ of $I(X)$ containing $i(X)$ for which 
$x, y, z \in E$ implies that $x \langle y | z \rangle \in E$.
In this case the $C^*$-algebra   ${\mathcal F}(X)$
in the last paragraph may be taken to be
the closed span in ${\mathcal D}$ of the 
set $\{ \langle y | z \rangle : y, z \in \partial X \}$.
Then $\partial X$ is a right $C^*$-module over ${\mathcal F}$.  
 
We end this preliminaries section with a Lemma:

\begin{lemma} \label{inc}  Let $Z$ be a right $C^*$-module
over a $C^*$-algebra $B$.  Then there is an injective
envelope $I(Z)$ of $Z$ which is a right $C^*$-module
over a $C^*$-algebra ${\mathcal R}$, with the following
properties: ${\mathcal R}$  contains  
$B$ as a $C^*$-subalgebra, and the 
module action $I(Z) \times {\mathcal R}$ restricted
to $Z \times B$ is the original one.
\end{lemma}

\begin{proof}    We may suppose that $B$ is unital.
Form the linking $C^*$-algebra ${\mathcal L}(Z)$ as 
in Section 2, and suppose that 
it is suitably represented nondegenerately as a 
$C^*$-subalgebra of $B(H \oplus K)$.  Thus
$1_B = I_K$.   The following is a mild variant 
of the Hamana-Ruan construction of 
the injective envelope, and the reader
may want to follow along with this  construction 
in any of the sources \cite{Ham2,RInj,BP,Pau}.
Consider the operator system $${\mathcal S}_B(Z) \; \; = 
 \; \; \left[ \begin{array}{ccl}  \Cdb I_H & Z \\ \bar{Z} & B
\end{array} \right]  \; \; \subset
B(H \oplus K). $$
Let $\Phi$ be a {\em minimal ${\mathcal S}_B(Z)$ projection}
on $B(H \oplus K)$, this is a completely positive idempotent map whose 
range is an injective envelope $I({\mathcal S}_B(Z))$
of ${\mathcal S}_B(Z)$ (see the cited references).
By a result of Choi and Effros \cite{CE},
 $I({\mathcal S}_B(Z))$ is a $C^*$-algebra with a
new product $x \circ y = \Phi(x y)$.   Also $I({\mathcal S}_B(Z))$ may
be regarded as a $2 \times 2$ matrix algebra with respect to
the canonical diagonal projections  $p = I_H \oplus 0$ and
$q = 0 \oplus I_K$.   Let ${\mathcal R}$ be the 
$2$-$2$-corner $q I({\mathcal S}_B(Z)) q$, this is
a $C^*$-subalgebra of $I({\mathcal S}_B(Z))$.  By definition 
of the new product it is easy to see that 
$B$ is a $C^*$-subalgebra of ${\mathcal R}$.  Let $E$ be the 
$1$-$2$-corner $p I({\mathcal S}_B(Z)) q$.  Clearly $E$
is injective, and is also a right 
$C^*$-module over ${\mathcal R}$.   In fact by definition
of the new product it is easy to see that the right
action of ${\mathcal R}$ on $E$ extends the action
of $B$ on $Z$.   By \cite[Theorem 2.6]{BP}, $E$ is injective in 
the category of operator $B$-modules.  We wish to show that 
$E$ is an injective envelope of $Z$.  To do this
we first show that $Id_E$ is the only completely contractive $B$-module
map $u : E \to E$ extending the identity map on
$Z$.  For  by Suen's variant on Paulsen's
lemma \cite{Sue}, such $u$ is the corner of a
completely positive map $\Psi$  on the 
subspace ${\mathcal S}_B(E)$ of $I({\mathcal S}_B(Z))$, 
such that $\Psi$ extends the identity map on ${\mathcal S}_B(Z)$. 
Extend $\Psi$  further to a complete contraction from $I({\mathcal S}_B(Z))$
to itself.  By the rigidity property of the
injective envelope, this latter map and hence
also $u$ must be the identity map.

The result is completed with an appeal to the fact
from  that the injective envelope of $Z$ is also the $B$-module 
injective envelope of $Z$
\cite[Theorem 2.6]{BP}.  The idea for this is
as follows: by that result in \cite{BP} 
 any injective envelope $I(Z)$ of $Z$ can 
be made into a $B$-module which is 
injective as an operator $B$-module.  A routine diagram chase,
 using facts from the last paragraph, shows 
that  $I(Z) \cong E$ completely isometrically and as $B$-modules.  \end{proof}

\section{Relative left multipliers between two spaces}
 
In this section $X$ and $Y$ are operator spaces possessing
noncommutative Shilov
boundaries $(\partial X,i)$ and $(\partial Y,j)$
respectively, which will be fixed for the
remainder of this section.  We will assume also
that $\partial X$ and $\partial Y$ are
right $C^*$-modules over $C^*$-algebras $B$
and $C$ respectively, where $B$ is $C^*$-subalgebra
of $C$.

The next result generalizes 
important facts
about left multipliers
on a single operator space.
To explain the notation in this result: the inner products
 in (ii) are the (right) Shilov inner products
on $Y$ and $X$ respectively, and  the matrices there
are indexed on rows by  $i$, and on columns by
$j$.   The first norm in (iii) is just the norm in $M_{2n,n}(Y)$,
the second is the norm on $M_n(X \oplus Y)$ 
inherited from $M_n(\partial X \oplus^C
\partial Y)$.
An explicit formula for this norm was given above Theorem
\ref{PasB2}.

Although the next result was stated in the survey \cite{Bsv}, 
it was proved only in a special case.
  
\begin{theorem} \label{chmu2}  Let
$X$, $Y, \partial X, \partial Y, B$ and $C$ be as above,
where $B \subset C$,
and let $T  \colon X \rightarrow Y$ be a linear map.
The following are equivalent:
\begin{enumerate}
\item [{\rm (i)}]   $T$   is the restriction
to $X$ of a (necessarily unique)
completely contractive right $B
$-module map $S :
\partial X \to \partial Y$.
\item [{\rm (ii)}]   $[\langle T(x_i) | T(x_j) \rangle] \leq [\langle x_i |x_j
\rangle ]$
for all $m \in \Ndb$ and $x_1, \cdots , x_m \in X$.
\item [{\rm (iii)}]  For all $n \in \Ndb$ and matrices
$[x_{ij}] \in M_n(X)$, $[y_{ij}] \in M_n(Y)$ we have
$$\left| \left|  \left[ \begin{array}{cl} T x_{ij} \\ y_{ij} \end{array}
\right] \right| \right|
\; \; \leq \; \;  \left| \left| \left[ \begin{array}{cl} x_{ij} \\ y_{ij}
\end{array}
\right] \right| \right|.$$
\end{enumerate}
If $B
= C$
then we may replace `completely contractive' by
`contractive' in {\rm (i)}.  \end{theorem}
 
\begin{proof}
(i) $\Rightarrow$ (ii) \ 
If $B = C$ then
we gave a simple proof of this implication in \cite{Bsv}.
In the general case we sketch another argument.
The point is that 
Paschke's proof of the implication (i) $\Rightarrow$ (ii) in
Theorem \ref{PasB2} has a matricial version
that  works here.  To cite details, we will use
the notation in Paschke's proof
(or in Section 10 of \cite{Bsv} where we reproduced
Paschke's proof).
We need to replace $h_n$ there by the matrix $H = 
([\langle x_i |x_j
\rangle] + \frac{1}{n} I_m)^{-\frac{1}{2}}$, $x_n$ by
the row matrix $v = [x_1 , \cdots , x_m] H$, and
the expression $\langle x | x \rangle$ by $[\langle x_i |x_j
\rangle]$.    One shows analogously to the proof 
in \cite{Bsv}
that $\Vert v \Vert \leq 1$.
Since $T$ is completely contractive, 
$T$ applied entrywise to $v$ has 
norm $\leq 1$; and then one proceeds 
along the earlier line.    We leave the details as an exercise. 

(ii) $\Rightarrow$ (iii) \  This is easy 
using the formula above Theorem \ref{PasB2} (see
e.g.\  \cite{Bsv}).
 
(iii) $\Rightarrow$ (i)   Note that
(iii) says that the map $T \oplus Id : X \oplus Y \to 
C_2(Y)$ is completely contractive, when $X \oplus Y$ is
viewed as a subspace of $\partial X \oplus^C
\partial Y$.   It follows
from a tedious diagram chase that any injective envelope
$I(\partial Y)$ of $\partial Y$ is an injective
envelope $I(Y)$ of $Y$.  Thus we may write
$I(Y)$ for the injective
envelope of $\partial Y$ in Lemma \ref{inc},
this is a $C^*$-module over ${\mathcal R}$,
where $C$ is a $C^*$-subalgebra of ${\mathcal R}$.
Then  $\partial X \oplus^C
\partial Y$ is
a subspace of 
$\partial X \oplus^{\mathcal R} I(Y)$. 
We may now follow the proof of the  implication (iii) $\Rightarrow$ (i)   
in  \cite{Bsv}, but with ${\mathcal D}$ replaced by 
${\mathcal R}$, to extend
$T$ to a right ${\mathcal R}$-module map
$\tilde{T} : I(X) \to I(Y)$.  Since $\tilde{T}$ is
also a right $B$-module map  we may conclude the 
proof as we did in \cite{Bsv}.
\end{proof}                      

For $X, Y$ as above
we define a {\em relative left multiplier} from $X$ to $Y$
to be a map
$T : X \to Y$ such that a positive scalar multiple of $T$
satisfies the equivalent conditions of Theorem
\ref{chmu2}.  We write $\M^{rel}_l(X,Y)$ for the set of such
relative left multipliers.
Note that $\M^{rel}_l(X,X)$ is simply the space $\M_l(X)$
in the Introduction.   
To define an operator space structure on 
$\M^{rel}_l(X,Y)$ we first observe that as in 
\cite[p.\ 303]{BShi} there is
a canonical linear isomorphism
$$\{ S \in CB_{B}
%{\mathcal F}(X)}
(\partial X,
\partial Y) : S(X) \subset Y \} \; \; 
%\longrightarrow 
\cong \; \;
 \M^{rel}_l(X,Y),$$
given by $S \mapsto S_{|X}$.  Since 
$CB(\partial X,\partial Y)$ is an operator space
so is the set on the left side of the 
last displayed expression.   We may therefore use
the linear isomorphism above to give $\M^{rel}_l(X,Y)$ an operator 
space structure.
Note that  Theorem \ref{chmu2} 
gives alternative descriptions of the unit ball of
$\M^{rel}_l(X,Y)$. 

In  \cite{Bsv} we gave some examples of relative left multipliers.
We also proved the following result, which we shall not
use in the present paper:

\begin{proposition} \label{mnml2}   If $X, Y$ are as in Theorem {\rm 
\ref{chmu2}},
and if $m,n \in \Ndb$  then we have
$M_{m,n}(\M^{rel}_l(X,Y)) \cong \M^{rel}_l(C_n(X),C_m(Y))$
 completely 
isometrically.
\end{proposition}

\section{Adjointable maps between two operator spaces}

In this section we consider two operator spaces $X, Y$  with fixed
noncommutative Shilov boundaries $\partial X$ and $\partial Y$
which are 
right $C^*$-modules over the same $C^*$-algebra $B$.

\begin{theorem} \label{genad}  Let $X, Y$ be as above
and suppose that $T : X \to Y$.  The following are equivalent:
\begin{itemize} 
\item [(i)]   $T$ is the restriction to $X$ of
an adjointable (in the usual $C^*$-module sense) 
$B$-module
 map $R : \partial X \to \partial Y$ such that $R(X) \subset Y$ and 
$R^*(Y) \subset X$,
\item [(ii)]  There exists a map $S : Y \to X$
such that $\langle T(x) | y \rangle = \langle x | S(y) \rangle$
(these are the (right)
Shilov inner products) for all $x, y \in X$.
\end{itemize}
Moreover the set $\A_l(X,Y)$ consisting of maps $T$ satisfying 
condition {\rm (ii)} above, is
a closed subspace of $B(X,Y)$ which is
a $C^*$-bimodule over the algebras
$\A_l(X)$ and $\A_l(Y)$.  The module actions here on
$\A_l(X,Y)$ are 
simply composition of operators.  The 
$\A_l(X)$-valued inner product 
on $\A_l(X,Y)$  is $\langle T | R \rangle = S R$, 
for $T, R \in \A_l(X,Y)$  
where $S$ is related to $T$ as in {\rm (ii)} above.
\end{theorem}                                
 
\begin{proof}   We leave it to the reader to check that
any $T \in \A_l(X,Y)$ is linear; that the 
map $S$ in (ii) is necessarily unique and linear;
that $\A_l(X,Y)$ is an $\A_l(Y)$-$\A_l(X)$-bimodule;
and  that the $\A_l(X)$-valued inner product specified above 
does indeed take values in $\A_l(X)$.  In fact the
only nontrivial part of the proof that 
$\A_l(X,Y)$ is a right $C^*$-module consists in showing
that for $T \in \A_l(X,Y)$,
(a) $T^* T \geq 0$ in $\A_l(X)$, and  
(b) $\Vert T^* T \Vert = \Vert T \Vert^2$.  Here 
$T^*$ denotes the map $S$ in (ii).  
In fact (a) follows from Theorem 4.10 (2) in \cite{BShi}, 
since 
$$\langle T^* T x | x \rangle = \langle T x | T x \rangle \geq 0,
\qquad x \in X .$$   To prove (b) we first note that if $R \in \A_l(X)_+$, 
with $R = V^* V$ for a $V \in \A_l(X)$, then
$$
%\Vert R \Vert =
\sup \{ \Vert \langle R x | x \rangle \Vert : 
%\Vert x \Vert \leq 1 \}
x \in {\rm Ball}(X) \} 
= \sup \{ \Vert \langle V x | V  x \rangle \Vert : 
%\Vert x \Vert \leq 1 \}
x \in {\rm Ball}(X) \} = \Vert V \Vert^2 = \Vert R \Vert.$$
Setting $R = T^* T$ and using (a) we see that
$\Vert T^* T \Vert$ equals 
$$\sup \{ \Vert \langle T^* T  x | x \rangle \Vert :
x \in {\rm Ball}(X) \}
= \sup \{ \Vert \langle T x | T  x \rangle \Vert :
x \in {\rm Ball}(X) \} = \Vert T \Vert^2.$$
%\Vert x \Vert \leq 1 \} 
It follows that   $\Vert T \Vert = \Vert T^* \Vert$
 as in the Hilbert space case.
%operators to
%verify also that $\Vert T \Vert = \Vert T^* \Vert$.

(i)  $\Rightarrow$ (ii)  This is obvious.

(ii) $\Rightarrow$  (i)  Suppose that $T$ satisfies (ii).
Then $T^* T \in \A_l(X)$ by the first part of the proof.
For $x_1, \cdots ,x_n \in X$ and $b_1, \cdots ,b_n \in B$,
we define $R(\sum_k x_k b_k) = \sum_k T(x_k) b_k$.  To see
that $R$ is well defined and bounded, set $u = \sum_k x_k b_k$,
take
$y_1, \cdots , y_m \in Y$ and $c_1, \cdots , c_m \in  B$
%z \in {\rm Ball}(\partial X)$ 
and set $v = \sum_k y_k c_k$.  Then
%compute:
\begin{equation} \label{adeq} 
\langle v | \sum_k T(x_k) b_k \rangle 
= \sum_{i,j} c_j^* \langle y_j | T(x_i) \rangle b_i
= \sum_{j} c_j^* \langle T^*(y_j) | u \rangle
= \langle \sum_k T^*(y_k) c_k | u \rangle .
 \end{equation}
%
%.$$

%= \sum_{j} f_j^* \langle T^* T(x_j) | \sum_i x_i b_i \rangle.$$
%= \sum_k \langle z | T(x_k) \rangle b_k 
%= \langle T^* z |\sum_k  x_k b_k \rangle.$$
Setting $y_k = T(x_k)$ and $c_k = b_k$ we obtain  from 
(\ref{adeq}) and a  Cauchy-Schwarz inequality:
$$\Vert 
%\langle
\sum_k T(x_k) b_k
%z | \sum_k T(x_k) f_k \rangle \Vert 
\Vert^2 \leq 
%\Vert T^* z \Vert
%\sum_k T(x_k) f_k \Vert^2 \leq 
\Vert \sum_k T^* T(x_k) b_k \Vert \; \Vert u \Vert .$$
% \Vert \sum_k x_k b_k  \Vert
%\leq \Vert T \Vert \Vert \sum_k x_k b_k  \Vert.$$
Now the mapping $\sum_k  x_k b_k  \mapsto \sum_k T^* T(x_k) b_k$
is simply the unique $B$-module map on $\partial X$ extending
$T^* T \in \A_l(X)$ (see \cite{BShi}
Theorem 4.10, in conjunction with 
the observation in equation (\ref{obre})),
 and this extension has the same norm.
Thus 
% that $T^* T $ extends uniqu 
$$\Vert \sum_k T(x_k) b_k    \Vert^2 \leq   \Vert T^* T \Vert \;
\Vert u \Vert^2 =  \Vert T \Vert^2 \; \Vert u \Vert^2.$$
%It follows that $$\Vert \sum_k T(x_k) b_k \Vert
%\leq \Vert T \Vert \Vert \sum_k x_k b_k  \Vert.$$
%%\leq  \Vert T^* T \Vert \Vert \sum_k x_k _k  \Vert^2.$$
Thus $R$ is bounded and well defined.

Since $R$ is bounded, it extends by density to a
unique bounded $B$-module
map $R : \partial X \to \partial Y$.
%, and it is easy to see that $R$ is a right ${\mathcal F}$-module map.
   Similarly $T^*$ extends to a   bounded map $S :
 \partial Y \to \partial X$.
We leave it as an exercise using (\ref{adeq}),
 to check that $R$ 
is adjointable with adjoint $S$, and satisfies (i).
\end{proof} 

As in the last proof we write $T^*$ for the unique $S$ related 
to $T$ in (ii) of the Theorem.  We call such maps $T$ {\em 
relatively adjointable}.      Strictly speaking
we should write $\A^{rel}_l(X,Y)$ for what we wrote as
$\A_l(X,Y)$ above, but for simplicity we use the
shorter notation in this section.
  As was the case for
$\M^{rel}_l(X,Y)$, the space $\A_l(X,Y)$ is only defined 
relative to fixed noncommutative Shilov boundaries $\partial X$ and
$\partial Y$.  There are frameworks in which
one may remove this
`relative' nature, as we shall see in the next sections. 
 
\medskip
 
{\bf Remarks: } 1)  It is easy to see that the set 
$$\{ R \in \Bdb_B(\partial X, \partial Y) : R(X) \subset Y , 
R^*(Y) \subset X \}$$
is a right $C^*$-module over the $C^*$-algebra
$$\{ R \in \Bdb(\partial X) : R(X) \subset X ,
R^*(X) \subset X \} .$$
By basic properties of 
`ternary morphisms' (see e.g.\ \cite{Ham2}),
the restriction map from this $C^*$-module 
onto $\A_l(X,Y)$ is a 
completely isometric surjective ternary isomorphism.

2)  If $W$ is a third operator operator space
whose noncommutative Shilov boundary $\partial W$ 
is also a right $C^*$-module over the same algebra
$B$
% ${\mathcal F}$
as above, then `composition  of operators' 
is a well defined bilinear 
map $\A_l(X,Y) \times \A_l(W,X) \to \A_l(W,Y)$.  Similar 
assertions hold for the $\M_l(\cdot,\cdot)$ spaces.
  
\section{Multipliers relative to a superspace}

In this section we consider a fairly general situation in which we
can remove some of
the relative nature of spaces $\M^{rel}_l(X,Y)$ and $\A^{rel}_l(X,Y)$
considered above. 
 
\begin{definition} \label{Shs}  Consider a
pair $(X,Y)$ of closed subspaces 
of an operator space $V$.   Suppose that
$X$ has the property that there is a noncommutative Shilov boundary $(\partial V,i)$ 
of $V$ such that
the smallest closed ${\mathcal F}(V)$-submodule of $\partial V$
containing $i(X)$ is a noncommutative Shilov boundary of $X$.
Suppose that $Y$ has the same property.
Then we say that $(X,Y)$ is  a 
{\em $\partial$-compatible $V$-pair}
\end{definition}
                                                   
By the universal property of the noncommutative Shilov boundary
\cite{Ham2,BShi}, together with \cite[Proposition 2.1 (iv)]{Ham2},
and a routine diagram chase, it is easy to see that the notions
above do not depend on the particular Shilov boundary of $V$ considered
above.   We will not use this, but one may rephrase the statement
``the smallest closed ${\mathcal F}(V)$-submodule of $\partial V$
containing $i(X)$ is a noncommutative Shilov boundary of $X$"
as a combination of two statements:
1) the `subTRO' $Z$ of $\partial V$ generated by 
$i(X)$ (see e.g.\ \cite{Ham2,Bsv} for the definition
of this) is a noncommutative Shilov boundary of $X$,
and 2) $Z$ is a right ${\mathcal F}(V)$-submodule of $\partial V$.
 
If $(X,Y)$ is a $\partial$-compatible $V$-pair, then
 %and that  ${\mathcal F}(V)$ is as defined in
we will henceforth in this section reserve the
symbols $\partial X$ and $\partial Y$
for the particular noncommutative Shilov boundaries
 of $X$ and $Y$ respectively mentioned in the last 
paragraph; these are submodules of $\partial V$.
%, and ${\mathcal F}(X)$ (or ${\mathcal F}$ when $X$ is
%understood) for the $C^*$-subalgebra of ${\mathcal F}(V)$
%generated by $i(X)^* i(X)$.  
It is easy to see
that $\partial X$ and $\partial Y$ are  right $C^*$-modules over
%${\mathcal F}(X)$, and also over 
${\mathcal F}(V)$.  Thus $\partial X \oplus_c
\partial Y$,  the $C^*$-module sum, is a
$C^*$-module over ${\mathcal F}(V)$.
We clearly have canonical complete isometric embeddings
$$X \oplus_V Y  \; \hookrightarrow \; \partial X \oplus_c
%{{\mathcal F}(Y)}
\partial Y  \; \hookrightarrow \; C_2(\partial V). $$   

The second matrix norm in (iii) below is the norm on $M_n(X \oplus_V Y)$.
 
\begin{corollary} \label{chmu3}  Let
$(X, Y)$ be a $\partial$-compatible $V$-pair, 
let $\partial X, \partial Y, {\mathcal F}(V)$ 
be as above, and set $C = {\mathcal F}(V)$.
If  $T  \colon X \rightarrow Y$ is a linear map
then the following are equivalent:
\begin{enumerate}
\item [{\rm (i)}]   $T$   is the restriction
to $X$ of a 
%completely check
contractive right $C$-module map $S :
 %{\mathcal F}(X) $-module map $S :
\partial X \to \partial Y$.
\item [{\rm (ii)}]   $[\langle T(x_i) | T(x_j) \rangle] \leq [\langle x_i |x_j
\rangle ]$
for all $m \in \Ndb$ and $x_1, \cdots , x_m \in X$.
\item [{\rm (iii)}]  For all $n \in \Ndb$ and matrices
$[x_{ij}] \in M_n(X)$, $[y_{ij}] \in M_n(Y)$ we have
$$\left| \left|  \left[ \begin{array}{cl} T x_{ij} \\ y_{ij} \end{array}
\right] \right| \right|
\; \; \leq \; \;  \left| \left| \left[ \begin{array}{cl} x_{ij} \\ y_{ij}
\end{array}
\right] \right| \right|.$$
\end{enumerate}
\end{corollary}

\begin{proof}   Follows from Theorem \ref{chmu2} 
with $B = C = {\mathcal F}(V)$. 
\end{proof}

\begin{definition} \label{Vmud}
If $(X, Y)$ is a $\partial$-compatible $V$-pair,
and if $T : X \to Y$ is such that a positive 
scalar multiple of $T$
satisfies the equivalent conditions of Corollary 
\ref{chmu3}, then we call $T$ a
{\em left $V$-multiplier} from $X$ to $Y$
We write $\M^V_l(X,Y)$, or $\M_l(X,Y)$ when $V$ is
understood, for the set of such
left $V$-multipliers.  
We write $\A_l^V(X,Y)$ for the $C^*$-bimodule in 
Theorem \ref{genad} (taking $B = {\mathcal F}(V)$ there).
The maps in $\A_l^V(X,Y)$ will be called left
{\em $V$-adjointable}.
\end{definition}

We will identify 
$\M^V_l(X,Y)$ with the operator space $\M^{rel}_l(X,Y)$
from Section 3, where 
the noncommutative Shilov boundaries of $X$ and $Y$
are taken to be the ones mentioned at the start of the current
Section.

\begin{proposition} \label{Vmcm}  Let $Y, Z$ be right $C^*$-modules
over a $C^*$-algebra $B$.  Set $V = Y \oplus_c Z$, and
regard $Y, Z$ as subspaces of $V$.   Then $(Y,Z)$ is a
$\partial$-comparable $V$-pair,
and $\M_l^V(Y,Z) \cong B_B(Y,Z)$ and $\A_l^V(Y,Z) \cong \Bdb(Y,Z)$.
\end{proposition}
 
\begin{proof}  Denote
 the closed span of the
 range of the canonical $B$-valued inner product on
$Y \oplus_c Z$ by ${\mathcal F}$.   In this case 
(see e.g.\ \cite{BShi}) one can take  
$\partial V = V$, viewed as a right $C^*$-module ${\mathcal F}$.
Then $Y, Z$ are also
$C^*$-modules over  ${\mathcal F}$, $(Y,Z)$ is a
$\partial$-comparable $V$-pair, and $\M_l^V(Y,Z) \cong B_{\mathcal F}(Y,Z)$ and
$\A_l^V(Y,Z) \cong \Bdb_{\mathcal F}(Y,Z)$.  We now may appeal
to the principle in equation (\ref{obre}).
\end{proof}

\begin{lemma}  \label{shbs}  If $(X,Y)$ is a $\partial$-compatible $V$-pair
then $\partial X \oplus_c \partial Y$
is a noncommutative Shilov boundary of $X \oplus_V Y$.
\end{lemma}  

\begin{proof}   First observe that $\partial X \oplus_c \partial Y$
is a left operator $\ell^\infty_2$-submodule of $C_2(\partial V)$.
The canonical map $X \oplus_V Y \to \partial X \oplus_c \partial Y
$ is a complete isometry as noted above.
Inside $\Bdb(\partial X \oplus_c \partial Y)$ there 
is a copy of $\ell^\infty_2$ (this is true for the 
sum of any two $C^*$-modules).
We may follow the proof in \cite[Theorem A.13]{BShi}:
one supposes that $W$ is a `ternary ideal' in $\partial X \oplus_c \partial Y$
such that the canonical map $X \oplus_V Y \to (\partial X \oplus_c \partial Y)/W$
is a  complete isometry, and then one needs to
show that $W = (0)$.   This is accomplished by letting $W_1 = e_1 W \subset 
\partial X$, and $W_2 = e_2 W \subset
\partial Y$, where $e_i$ is the `standard basis' for $\ell^\infty_2$, and 
showing that the canonical maps $X \to 
(\partial X)/W_1$ and $Y \to
(\partial X)/W_2$ are complete isometries.  Since the reasoning is identical
to that in \cite[Theorem A.13]{BShi} we omit the details.  
\end{proof}  

%{\bf Remark: }
We write $\epsilon_X$ and $P_X$ for the canonical
inclusion and projection maps between $X$ and $X \oplus_V Y$.
Similarly for $\epsilon_Y$ and $P_Y$.  These maps are
restrictions of the canonical adjointable
inclusion and projection maps between the 
$C^*$-module $\partial X \oplus_c \partial Y$ and its summands.
It is clear from the definitions in  \cite{BEZ} that 
$p = \epsilon_X \circ P_X$ is a left $M$-projection on
$X \oplus_V Y$, onto the right $M$-summand $X \oplus 0$.
% (see \cite{BEZ} for definitions of these terms).
  Similarly $q = \epsilon_Y \circ P_Y$ is the left $M$-projection onto
$0 \oplus Y$.

For $X, Y, V$ as above, $X \; \oplus \; \{ 0 \}$ and $\{ 0 \}
 \;  \oplus \; Y$
 are a $\partial$-compatible $X \oplus_V Y$-pair, as may be seen using
Lemma \ref{shbs}.
Thus it seems that in most situations we may assume without
loss of generality that $X, Y$ are complementary right $M$-summands
in $V$ (by `replacing' $V$ by $X \oplus_V Y$, and using the
observation in the last paragraph).

Conversely, if $(X,Y)$ is a
$\partial$-compatible $V$-pair, and if also $X$ and $Y$ are 
`complementary' right $M$-summands in $V$, then $V \cong
X \oplus_V Y$ completely isometrically.  This is because,
by one of the definitions of a right $M$-summand,
the map $$V \; \to C_2(V)  \; : \;  x + y \; \mapsto \; \left[
\begin{array}{c}
x \\
y \end{array}
\right]$$
is a complete isometry, and its range 
is the space $X \oplus_V Y$ defined in the Introduction.

The next observation we make is that since the operator algebra
$\M_l(X \oplus_V Y)$ (resp. $C^*$-algebra $\A_l(X \oplus_V Y)$)
contains the two canonical complementary  
projections $p, q$ mentioned a few paragraphs
above, it splits as a $2 \times 2$ matrix algebra
(resp. $C^*$-algebra).  We first claim that 
the $1$-$1$-corner is completely isometrically
homomorphic to $\M_l(X)$ (resp. $\A_l(X)$). 
 To see this 
%we will use the simple  .  
consider the map $\theta : \M_l(X) \to \M_l(X \oplus_V Y)$ taking $T$ to
$\epsilon_X \circ T \circ P_X$.  This may be viewed as the 
restriction to $\M_l(X)$ of the map $R \to 
\epsilon_{\partial X} \circ R \circ P_{\partial X}$ from 
the space of bounded module maps on $X$,
to the space of bounded module maps on $X \oplus_c \partial Y$.
Thus it is a well defined completely contractive homomorphism.
From this argument, or directly,
it is easy to see that $\theta$ is completely isometric.
If $S = p S' p$ for a map $S' \in \M_l(X \oplus_V Y)$,
then by the Lemma \ref{shbs} and Theorem 
\ref{chmu2}, $S'$
%_{\vert X}$ 
is the restriction to $X \oplus_V Y$ of 
%the composition of three 
a bounded
${\mathcal F}(V)$-module map $R'$ on $\partial X \oplus_c \partial Y$.
Then $S$ is the restriction to $X \oplus_V Y$ of 
$\epsilon_{\partial X} \circ P_{\partial X} \circ R' \circ \epsilon_{\partial X} \circ 
P_{\partial X}$.  From this it is clear that
$\theta(P_X \circ S' \circ \epsilon_X) = S$.
%and $\theta(S_{\vert X}) = S$.  
Thus $\theta(\M_l(X)) = p \M_l(X \oplus_V Y) p$.
If $R \in \Bdb(X)$ then $\epsilon_{\partial X} \circ R \circ P_{\partial X}
\in \Bdb_{\cdot}(\partial X \oplus_c \partial Y)$.  Thus it is
easy to argue that $\theta$ induces a $*$-monomorphism from $\A_l(X)$ onto 
$p \A_l(X \oplus_V Y)p$.

By identical reasoning we have completely isometries from 
$\M_l(Y)$, $\M_l(X,Y)$ and $\M_l(Y,X)$ into the other 
three corners of $\M_l(X \oplus_V Y)$.
Similar assertions hold for the $\A_l(\cdot)$ spaces.

\medskip
 
The following is the analogue of another important 
property of left $V$-multipliers on a single space
\cite{BShi}.  It can be stated in
many forms, but perhaps the following is the most
concise:
 
\begin{proposition} \label{chmuiv}  If  $X, Y$
are complementary right $M$-summands of an
operator space $V$ (see the discussion after Lemma
\ref{shbs}),
and if $(X,Y)$ is a $\partial$-compatible $V$-pair, then a
linear map $T : X \to Y$ is a left $V$-multiplier
if and only if there is a completely
isometric linear embedding of $V$ into
a $C^*$-algebra $A$,
and an $a \in {\rm Ball}(A)$ with $T x = a x$ for all $x \in X$.
\end{proposition}
 
\begin{proof}
 ($\Rightarrow$) \
Suppose that $T : X \to Y$ is a left $V$-multiplier.
Since $\M_l^V(X,Y)$ may be regarded as a corner
of $\M_l(X \oplus_V Y)$, and since $X \oplus_V Y
\cong V$ by the discussion after Lemma
\ref{shbs}, we may regard $T$ as a left multiplier 
$R$ of $V$.  Thus by the `one-space variant' of the
result we are trying to prove,
%\ref{chmu}
there is a completely
isometric linear embedding $\sigma : V \to A$,
and an $a \in {\rm Ball}(A)$ with $\sigma(R(x)) = a \sigma(x)$,
for all $x \in X$.  However $R(x) = T(x)$.

($\Leftarrow$) \  It is easy to show that 
the condition here implies Theorem \ref{chmu3} (iii).
  \end{proof}

We recall from \cite{Dual,BEZ}
that $\A_l(X)$ is
%are respectively a dual operator algebra and 
a $W^*$-algebra if $X$ is a dual operator space.  The
following generalizes this important fact: 

\begin{corollary} \label{iswm}  If $(X, Y)$ is
a $\partial$-compatible $V$-pair,
where $V$ is a dual operator space and 
$X, Y$ are weak* closed  subspaces of $V$, then 
$\A_l^V(X,Y)$ is a $W^*$-module.  Moreover, every
$T \in \A_l^V(X,Y)$ is automatically weak* continuous.
\end{corollary}

\begin{proof} 
To see that $\A_l^V(X,Y)$ is a $W^*$-module it suffices to show
that $\A_l^V(X,Y)$ is a dual space.  However
as we just saw, $\A_l^V(X,Y)$ is a `corner' in 
$\A_l(X \oplus_V Y)$.   Thus by  the fact
 mentioned above the Corollary,   
it suffices to show that $X \oplus_V Y$ is a dual operator
space.  However this is clear since $C_2(V)$ is
a dual operator space, and $X \oplus_V Y$ is 
easily seen to be weak* closed in $C_2(V)$.

The last assertion follows from the 
analogous fact for $\A_l(X \oplus_V Y)$ 
(see \cite{BEZ}), together with the
fact that the canonical inclusion and projection maps
between $X \oplus_V Y$ and its summands are 
 weak* continuous in this case (which follows from
basic operator space theory).   
\end{proof}  

The last result should be useful in the way that its
 `one-space predecessor' was (see e.g.\ \cite{BZ}).
  For example structural 
properties in a $W^*$-module (for example 
those considered in \cite{Har} or \cite{Rtro}) should have
implications for the pair $X, Y$.

It is often useful that the adjointable maps on 
an operator space,  or even on a Hilbert space,
are characterizable as the span of the Hermitian
(i.e.\ self-adjoint) ones.  The following may be viewed as the
`two-space' analogue of this fact.   
 
\begin{corollary} \label{lam3}  Let  $(X,Y)$ be a
$\partial$-compatible $V$-pair, set $B = {\mathcal F}(V)$,
and let $T : X \to Y$ be a linear map.  Then $T$ 
satisfies the equivalent conditions in Theorem 
\ref{genad} if and only if there is a  map $S : Y \to X$
such that the map  $(x,y) \mapsto (S(y),T(x))$ is a 
Hermitian in the Banach algebra $\M_l(X \oplus_V Y)$.
In this case $T^* = S$.
\end{corollary}
 
\begin{proof}   We leave this as an exercise.  The idea
is very similar to the last proof and the
discussion above it.  That is, using the canonical 
inclusion and projection maps, we transfer the 
desired statement to a statement about maps between
the $C^*$-modules
$\partial X, \partial Y$    and $\partial X \oplus_c \partial Y$.
\end{proof}

\section{Multipliers and  the injective envelope}  

In this brief section we list some variants of results 
in the last section, but with the noncommutative Shilov 
boundary replaced by the injective envelope.
 
We consider an operator space $V$, and 
fix an injective envelope $(I(V),i)$ of $V$, which is a right 
$C^*$-module over a $C^*$-algebra 
${\mathcal D} = {\mathcal D}(V)$ 
(see the paragraph before Lemma \ref{inc}).
  We say that a subspace $X$ of 
$V$ is a (right) {\em ${\mathcal D}$-subspace} if
there is a ${\mathcal D}$-submodule $W$ of $I(V)$ such that
$(W,i)$ is an injective envelope of $X$.  We call 
a pair $(X,Y)$ of  ${\mathcal D}$-subspaces of $V$,
an {\em $I$-compatible $V$-pair}.  
Clearly  $W$ is 
a right $C^*$-module over ${\mathcal D}(V)$ too.
We will write $W$ as $I(X)$, and ${\mathcal D}(X)$
for the $C^*$-subalgebra $W^* W$ of ${\mathcal D}(V)$.   Similar
notations hold for $Y$.  
In (i) below the notation $I_{11}(V)$ is used
 precisely in the sense of \cite{BP}.       
  
\begin{theorem} \label{newest} \cite{Bsv} Let $(X, Y)$ be an $I$-comparable
 $V$-pair.
Suppose further
that ${\mathcal D}(Y) \subset {\mathcal D}(X)$,
where  ${\mathcal D}(\cdot)$ is as defined
above.
%, {\mathcal D}$ be as above.  
If $T : X \to Y$ is a linear map,
then 
the following are equivalent:
\begin{itemize} \item [(i)]  There is an $a \in {\rm Ball}(I_{11}(V))$
such that $i(T x) = a i(x)$ for all $x \in X$.
\item [(ii)]  $T \oplus Id_X : C_2(X) \to Y \oplus_V X$ is
completely contractive.
\item [(iii)]  $T \oplus Id_Y : X \oplus_V Y \to C_2(Y)$  is
completely contractive.
\item [(iv)]   There is a $C^*$-algebra $A$,  a
completely isometric embedding $V \hookrightarrow A$, and
an $a \in  {\rm Ball}(A)$ such that $T x = a x$ for all $x \in X$.
\item [(v)] $T$ is the restriction to $X$ of a
contractive ${\mathcal D}(V)$-module map $S : I(X) \to I(Y)$.
\end{itemize}     
 
If further, $(X,Y)$ is a {\em $\partial$-compatible $V$-pair},
then the equivalent conditions above are also
equivalent to conditions {\rm (i)}--{\rm (iii)} in
Theorem \ref{chmu3}.
\end{theorem}

We next claim that
the discussion in the paragraphs between Lemma \ref{shbs} 
and Proposition \ref{chmuiv} above,
 is also valid for $I$-compatible $V$-pairs.
To see this one needs the following result:

\begin{lemma}  \label{Ibs}  If $(X,Y)$ is an $I$-compatible $V$-pair
then $I(X) \oplus_c I(Y)$
is an injective envelope for $X \oplus_V Y$.
\end{lemma}
 
\begin{proof}   Since $I(X)$ is a ${\mathcal D}(V)$-submodule 
of $I(V)$ it is a right $M$-ideal in $I(V)$ by 
\cite[Theorem 6.6]{BEZ}.  Since $I(X)$ is injective there is
a contractive projection $I(V)$ onto $I(X)$.
It follows from Theorem 3.10 (c) and Theorem 6.6
in \cite{BEZ} that 
%$I(X)$ is a right $M$-summand in $I(V)$ and indeed
there exists a contractive ${\mathcal D}(V)$-module map
projection from $I(V)$ onto $I(X)$.
Similarly for $I(Y)$.  It follows from 
Theorem \ref{PasB2} (iii) for example,
that there is
a contractive ${\mathcal D}(V)$-module map
projection from $C_2(I(V))$ onto $I(X) \oplus_c I(Y)$.
Thus $I(X) \oplus_c I(Y)$ is 
injective, since $C_2(I(V))$ is injective.
It suffices, by one of the equivalent definitions
of the injective envelope \cite{Ham2,RInj}, to show that
if $P$ is a completely contractive projection 
on $I(X) \oplus_c I(Y)$ which restricts to the
identity on $X \oplus_V Y$, then $P$ is the identity map.
If $\epsilon_X, P_X$ are as in the discussion below
 Lemma \ref{shbs}, then $P_{I(X)} \circ P \circ \epsilon_{I(X)}$
 is a complete contraction on $I(X)$ 
which restricts to $Id_X$.  By rigidity (see Section 2), 
$P_{I(X)} \circ P \circ \epsilon_{I(X)} = Id_{I(X)}$.
Similarly $P_{I(Y)} \circ P \circ \epsilon_{I(Y)}
= Id_{I(Y)}$.  Since $P^2 = P$, by pure algebra we must 
conclude that $P_{I(Y)} \circ P \circ \epsilon_{I(X)}$
and $P_{I(X)} \circ P \circ \epsilon_{I(Y)}$ are zero.  
Thus $P = Id$.
\end{proof} 

Since (by the Lemma)
the discussion in the paragraphs after Lemma \ref{shbs}
transfers to the present setting,
one may check that the conclusions of  Corollary \ref{iswm} are
true for $I$-compatible $V$-pairs too.

%{\bf Remark:} 
%We end this section with 
There is another characterization 
of left $V$-multipliers which is also analogous to the 
formulation of left multipliers in \cite{BP}.  
To state this characterization we suppose 
%as in the previous theorem
that $(X, Y)$ is an $I$-compatible $V$-pair.
For simplicity we also suppose that
${\mathcal D}(Y) \subset {\mathcal D}(X)$.
We then have as above 
that $I(X)$ and $I(Y)$ are 
right $C^*$-modules over ${\mathcal D}(X)$, and 
hence also over the $C^*$-algebra multiplier algebra
$M({\mathcal D}(X))$.  The 
latter $C^*$-algebra is injective too, by \cite[Corollary 1.8]{BP}.
Indeed $M({\mathcal D}(X)) \cong I_{22}(X)$ in the language of \cite{BP};
and henceforth we shall just write $I_{22}$ for $M({\mathcal D}(X))$.
We consider the `generalized linking $C^*$-algebra'
$A = \Bdb_{I_{22}}(I(Y) \oplus_c I(X) \oplus_c I_{22})$.
With respect to the canonical diagonal projections corresponding
to the identities of $\Bdb(I(Y)), \Bdb(I(X))$ and $I_{22}$ respectively,
 $A$ may be written as a $3 \times 3$ matrix $C^*$-algebra,
whose $k$-$\ell$-corner we write as $I_{k \ell}$, for $k, \ell \in
\{ 0,1,2 \}$.   Clearly $I_{02} = I(Y)$ and $I_{12} =
I(X)$.    With a little work 
one can show that $A$, and consequently also
$I_{ij}$, is injective.  We will not use this here however.
   We write $i$ and $j$ for the canonical maps
from $Y$ and $X$ into $I_{02}$ and $I_{12}$ respectively.
  
\begin{theorem} \label{i01}  Suppose that $(X, Y)$ is as
in the first part of
Theorem \ref{newest}.  Then a linear map $T : X \to Y$
satisfies conditions {\rm (i)}--{\rm (v)} in that 
Theorem if and only if there exists an element 
$a \in I_{01}$ such that $i(T x) = a j(x)$ for all 
$x \in X$.    
\end{theorem}
     
\begin{proof}   By \cite[Corollary 2.7 (iii)]{BP}
we have that 
$$I_{01} \cong \Bdb_{I_{22}}(I(X),I(Y)) = B_{I_{22}}(I(X),I(Y)).$$  Hence 
$$I_{01} \cong B_{I_{22}}(I(X),I(Y)) = B_{{\mathcal D}(X)}(I(X),I(Y))
= B_{{\mathcal D}(V)}(I(X),I(Y)),$$ 
using the principle in equation (\ref{obre}).   The result is clear from this
and Theorem \ref{chmu3} (v).  \end{proof}

This result, and the matching part of the last theorem,
may also be proved by a variation
of the proof given in \cite{Pau} of the analogous assertion for
$\M_l(X)$.

It should be interesting and useful
to extend other known results about $\M_l(X)$ and $\A_l(X)$
(for example those in \cite{BZ})
to the case of two spaces $X$ and $Y$.
         
\medskip
 
{\bf Acknowledgments:}  I am grateful to
Vrej Zarikian for comments and suggestions.

 \end{document}